\documentclass[10pt]{article}
\usepackage{amssymb}
\usepackage{hyperref}
\usepackage{latexsym}
\usepackage{color}
\usepackage{amsmath}
\usepackage{graphicx}
\usepackage{amsfonts}
\usepackage{amsthm}
\usepackage{graphicx}

\newtheorem{theorem}{Theorem}

\newtheorem{definition}{Definition}
\newtheorem{proposition}{Proposition}

\newtheorem{remark}{Remark}

\everymath{\displaystyle}
\begin{document}

\title{ On the generalized oblate spheroidal wave functions and applications}
\author{Tahar Moumni$^{a}$\thanks{%
Corresponding author: E-mail: moumni.tahar1@gmail.com}  and Ammari Amara$^{b}$  \\
$^{a}$ University of Carthage, Faculty of Sciences of Bizerte,\\
Department of Mathematics, Jarzouna,
Tunisia.\\
$^{b}$ University of Manar, Faculty of sciences of Tunis, \\
Department of Mathematics, Tunis, Tunisia.\\
} \maketitle

\noindent\textbf{Abstract} In this paper, we introduce a new set of
functions, which have the property of the completeness over a finite
and infinite intervals. This family of functions, denoted for
simplicity GOSWFs, are a generalization of the oblate spheroidal
wave functions. They generalize also the Jacobi polynomials in some
sens. The GOSWFs are nothing but the eigenfunctions of the finite
weighted bilateral Laplace transform
$\mathcal{F}_c^{(\alpha,\beta)}.$ We compute this functions by two
methods:  In the first one we use  a differential operator
$\mathcal{D}$ which commutes with $\mathcal{F}_c^{(\alpha,\beta)}.$
In the second one we use  the Gaussian quadrature method. As an
application, we use the GOSWFs to approximate bilateral weighted
Laplace  bandlimited functions and we show that they are more
advantageous then other classical basis of
$L^2((-1,1),(1-x)^\alpha(1+x)^\beta)dx$. Finally, we provide the
reader by some numerical examples that illustrate the theoretical
results.

\vskip0.5cm

\begin{quotation}
\noindent

\end{quotation}

\section{Introduction}
Prolate spheroidal wave functions (PSWFs), were first known as the
eigenfunctions of the following differential operator, see
\cite{Niven},
\begin{equation}\label{diff}
L_c(\varphi)(x)=(1-x^2)\frac{d^2\varphi}{dx^2}-2x\frac{d\varphi}{dx}-c^2x^2\varphi
.
\end{equation}
80 years later, D. Slepian, H. Landau and H. Pollack, see
\cite{Slepian-Pollak,landau pollak2, Landpoll} have shown that the
previous differential operator commutes with the following integral
operator
\begin{equation}\label{sinc}
F_c(f)(x)=\int_{-1}^{1}\frac{\sin c(x-y)}{\pi (x-y)}f(y)dy.
\end{equation}
Also, they  have given many properties of the PSWFs of high
practical importance, notably the double orthogonality over finite
 and infinite interval.
 Three years later, D.  Slepian showed in  \cite{slepian64}  that higher
dimensional construction of the PSWFs is  reduced to the $2D-$case.
Hence, the study of the $2D-$case is potentially important, see
\cite{Yoel}. He showed that the 
circular prolate spheroidal wave function (CPSWFs, for short), is an
eigenfunction of the finite Hankel transform.  That is
\begin{equation}\label{0.1}
\int_0^{1}\sqrt{cxy}J_\nu(cxy)\psi_{n,c}^\nu(y)dy=\gamma_{n,\nu}(c)\psi_{n,c}^\nu(x),\qquad
\nu \in\mathbb{N},\quad x\in(0,+\infty).
\end{equation}
 For more details about these functions and their computational methods, the reader is referred to
 \cite{slepian64}. As it is
done for the case of  the PSWFs, Slepian has proved that the
following differential operator
\begin{equation}\label{diff}
\mathcal{L}_c(\varphi)(x)=(1-x^2)\frac{d^2\varphi}{dx^2}-2x\frac{d\varphi}{dx}+(-c^2x^2+\frac{\nu^2-1/4}{x^2})\varphi
\end{equation}
commutes with the finite Hankel transform. Note here, that the
differential operators $\mathcal{L}_c$ is  a perturbation of the
differential operator
\begin{equation}\label{diff}
\mathcal{L}_0(\varphi)(x)=(1-x^2)\frac{d^2\varphi}{dx^2}-2x\frac{d\varphi}{dx}+\frac{\nu^2-1/4}{x^2}\varphi
\end{equation}
whose eigenfunctions are nothing but
\begin{equation}\label{tknu}
T_{k,\nu}(x)=\sqrt{2(2k+\nu+1)}x^{\nu+\frac{1}{2}}P_{k}^{(\nu,0)}(1-2x^2),\:\:k\geq
0\end{equation}

 where for $\alpha, \beta >-1,$ $P_{n}^{(\alpha ,\beta)}(x)$  is
the Jacobi polynomial
 of degree $n.$  In
\cite{karoui-moumni2}, we have extended the CPSWFs to an arbitrary
order $\nu>-1$ instead of $\nu\in\mathbb{N}$ and we showed that
CPSWFs share the same property as the PSWFs notably the double
orthogonality over finite and infinite interval.

Recently, in \cite{Wang1}, the authors have given a new extension of
the PSWFs that generalizes the Gegenbauer polynomials to an
orthogonal system with an intrinsic tuning parameter $c > 0.$ These
generalized PSWFs denoted by GPSWFs are defined as the
eigenfunctions of a Sturm-Liouville problem ${\mathcal D_x},$ that
commutes with   an integral operator ${\mathcal F_c^{(\alpha)}}.$
These operators are defined as follows,
\begin{eqnarray*}
\mathcal D_x u&=& -(1-x^2)^{-\alpha} \frac{d}{dx}
\left((1-x^2)^{\alpha+1}\frac{du}{ dx} \right)+(c^2 x^2)u,\quad
\alpha>-1,\quad c>0,\quad x\in
(-1,1),\\
\mathcal F_c^{(\alpha)}(\varphi)(x)&=&\int_{-1}^1 e^{icxt}
\varphi(t) (1-t^2)^{\alpha} \, dt.
\end{eqnarray*}
The authors have shown that the GPSWFs share similar properties with
the PSWFs. Also, they have presented
 a number of analytic and asymptotic formulae for the GPSWFs as well as their associated eigenvalues. They
  introduced
efficient algorithms for their evaluations.

Remark here that in the special case where $c=0$ the PSWFs are
nothing else the Legendre polynomials, the kth CPSWFs coincides with
$T_{k,\nu}(x)$ given by (\ref{tknu}) and the kth GPSWFs of Wang is
the kth Gegenbauer polynomial. In the literature, the PSWFs and all
its above mentioned generalization are called Slepian's functions,
this is perhaps due to the "lucky accident" of David Slepian.

A natural question that can be asked is the following: If there is a
kind of Slepian's functions that generalize the Jacobi polynomials.
In this paper, we answer positively to the previous question and we
show that the bounded solutions of the following equation

\begin{equation}  \label{17.0100}
(1-x^2)\frac{\partial^2\psi_{n}^{(\alpha,\beta)}(x; c)}{\partial^2
x}+(\beta-\alpha-(\alpha+\beta+2)x)\frac{\partial\psi_{n}^{(\alpha,
\beta)}(x; c)}{\partial  x} \notag +(c^2x^2
+cx(\beta-\alpha))\psi_{n}^{(\alpha,\beta)}(x;
c)=-\chi_{n}^{(\alpha,\beta)}(c)\psi_{n}^{(\alpha,\beta)}(x; c).
\end{equation}
generalizes in some  sense  some of the previous mentioned Slepian's
functions. In fact, the GOSWFs coincides with
\begin{itemize}
\item the PSWFs if $\alpha=\beta=0$ and $c= i\widetilde{c}$ and therefore with Legendre
polynomials if $\alpha=\beta=c=0.$ Here $i^2=-1$ and
$\widetilde{c}>0.$
\item Wang's GPSWFs if $\alpha=\beta$ and $c=i\widetilde{c}$ and therefore with
Gegenbauer polynomials if  $\alpha=\beta$ and $c=0.$
\item Jacobi polynomials if $c=0.$
\item  Chebyshev polynomials of first kind if
$\alpha=\beta=\frac{-1}{2}$ and $c=0.$
\item  Chebyshev polynomials of second kind if
$\alpha=\beta=\frac{1}{2}$ and $c=0.$
\end{itemize}

As an application of the GOSWFs, we show that the GOSWFs are more
advantageous in approximate the bilateral weighted Laplace
bandlimited functions  then other classical basis of
$L^2((-1,1),(1-x)^\alpha(1+x)^\beta dx)$ such as Jacobi polynomials.
We give also several numerical results that illustrate the
theoretical ones for different values of the parameters $\alpha,$
$\beta$ and $c.$

\noindent The outline of the paper is as follows: In section 2 we
recall some mathematical preliminaries about Jacobi polynomials and
Whittaker functions which will be used frequently later. Section 3
is devoted to define the GOSWFs and present some of their
properties. The goal of section 4 is the computation of these
functions and their corresponding eigenvalues by two different
methods. Finally, in section 5, we use the GOSWFs  to approximate
bilateral weighted Laplace bandlimited functions and we give some
numerical results.

\section{Mathematical preliminaries}

\subsection{Jacobi polynomials}

It is well known that for any two real numbers $\alpha ,\beta >-1,$
the Jacobi polynomial $P_{n}^{(\alpha ,\beta )}(x)$ of degree $n$ is
given by the following Rodriguez formula, see \cite{shen-tang}
\begin{equation}
P_{n}^{(\alpha ,\beta )}(x)=\frac{(-1)^{n}}{2^{n}n!}(1-x)^{-\alpha
}(1+x)^{-\beta }\frac{d^{n}}{dx^{n}}\left( (1-x)^{\alpha
+n}(1+x)^{\beta +n}\right) ,n\geq 0.  \label{eq2.1}
\end{equation}%
Let $k_{n}^{\alpha ,\beta }$ be the leading coefficient of   $%
P_{n}^{\alpha ,\beta }(x)$. Then $k_{n}^{\alpha ,\beta
}=\frac{\Gamma (2n+\alpha +\beta +1)}{2^{n}\cdot n!\Gamma (n+\alpha
+\beta +1)}.$ Moreover, if
\begin{equation*}
a_{n}^{\alpha ,\beta }=\frac{2^{(\alpha +\beta +1)}\Gamma (n+\alpha +1)\Gamma (n+\beta +1)}{%
n!(\alpha +\beta +2n+1)\Gamma (\alpha +\beta +n+1)},\:\:n\geq 0,
\end{equation*}%
then the set $\{\mathbf{P}_{n}^{(\alpha ,\beta )}(x),\;n\in
\mathbf{N}\}$ is an orthonormal basis of
$(L^{2}(-1,1),d\omega_{\alpha,\beta} )$, where
$d\omega_{\alpha,\beta} (x)=(1-x)^{\alpha }(1+x)^{\beta }dx$ and
$\mathbf{P}_{n}^{(\alpha ,\beta )}(x)=\frac{1}{\sqrt{a_{n}^{\alpha
,\beta }}}P_{n}^{(\alpha ,\beta )}(x)$. Note here that the
normalized Jacobi polynomial $\mathbf{P}_{k}^{(\alpha ,\beta )}(x)$
satisfy the following differential equation
\begin{eqnarray}
\mathcal{D}_{x}\mathbf{P}_{k}^{(\alpha ,\beta )}(x)
&=&(1-x^{2})\partial_x^2\mathbf{P}_{k}^{(\alpha ,\beta )}(x) +(\beta
-\alpha -(\alpha +\beta +2)x)\partial_x\mathbf{P}_{k}^{(\alpha
,\beta )}(x)  \notag  \label{17.0100} \\
&=&-\chi _{k}^{(\alpha ,\beta )}(0)\mathbf{P}_{k}^{(\alpha ,\beta
)}(x),
\end{eqnarray}%
where $\chi _{k}^{(\alpha ,\beta )}(0)=k(k+\alpha +\beta +1).$ To
proceed further, we recall the recurrence relation satisfied by
$\mathbf{P}_{k}^{(\alpha ,\beta )}(x),$

\begin{equation}
\mathbf{P}_{n+1}^{(\alpha ,\beta
)}(x)=(A_{n}x-B_{n})\mathbf{P}_{n}^{(\alpha ,\beta
)}(x)-C_{n}\mathbf{P}_{n-1}^{(\alpha ,\beta )}(x),  \label{Jacrecu}
\end{equation}%
where
\begin{eqnarray}
A_{n} &=&\sqrt{\frac{a_{n}^{\alpha ,\beta }}{a_{n+1}^{\alpha ,\beta }}}\frac{%
(2n+\alpha +\beta +1)(2n+\alpha +\beta +2)}{2(n+1)(n+\alpha +\beta
+1)}
\label{1***} \\
B_{n} &=&\sqrt{\frac{a_{n}^{\alpha ,\beta }}{a_{n+1}^{\alpha ,\beta }}}\frac{%
(\beta ^{2}-\alpha ^{2})(2n+\alpha +\beta +1)}{2(n+1)(n+\alpha
+\beta
+1)(2n+\alpha +\beta )}  \label{2***} \\
C_{n} &=&\sqrt{\frac{a_{n-1}^{\alpha ,\beta }}{a_{n+1}^{\alpha ,\beta }}}%
\frac{(\alpha +n)(\beta +n)(2n+\alpha +\beta +2)}{(n+1)(n+\alpha
+\beta +1)(2n+\alpha +\beta )}.  \nonumber
\end{eqnarray}

 Note here that  (\ref{Jacrecu}) can be written as follows
\begin{equation}
x\mathbf{P}_{k}^{(\alpha ,\beta )}=\alpha
_{k}\mathbf{P}_{k+1}^{(\alpha ,\beta )}+\beta
_{k}\mathbf{P}_{k}^{(\alpha ,\beta )}+\gamma
_{k}\mathbf{P}_{k-1}^{(\alpha ,\beta )}. \label{22}
\end{equation}%
where
\begin{equation*}
\alpha _{k} =\frac{1}{%
A_k}, \:\: \beta _{k} =\frac{B_k}{A_k},\:\: \gamma _{k}
=\frac{C_k}{A_k}.
\end{equation*}%

By multiplying both sides of (\ref{22}) by $x$ and using (\ref{22}), one gets
\begin{eqnarray}
x^{2}\mathbf{P}_{k}^{(\alpha ,\beta )} &=&\alpha _{k}\alpha
_{k+1}\mathbf{P}_{k+2}^{(\alpha ,\beta )}+\alpha _{k}(\beta
_{k+1}+\beta _{k})\mathbf{P}_{k+1}^{(\alpha ,\beta )}+(\alpha
_{k}\gamma _{k+1}+(\beta _{k})^{2}+\gamma _{k}\alpha
_{k-1})\mathbf{P}_{k}^{(\alpha ,\beta )}  \notag \\
&&+\gamma _{k}(\beta _{k}+\beta _{k-1})\mathbf{P}_{k-1}^{(\alpha
,\beta )}+\gamma _{k}\gamma _{k-1}\mathbf{P}_{k-2}^{(\alpha ,\beta
)},\quad \forall k\geq 2. \label{23}
\end{eqnarray}

\subsection{Whittaker Functions}

The Whittaker functions see \cite{Nilst} arise as solutions to the
Whittaker differential
equation.%
\begin{equation}\label{*8*}
\displaystyle\frac{d^{2}u}{dz^{2}}+(-\frac{1}{4}+\frac{\lambda }{z}+\frac{%
\frac{1}{4}-\mu ^{2}}{z^{2}})u=0.
\end{equation}%
Two solutions are given by the Whittaker functions $M_{\lambda
,\mu}(z)$, $W_{\lambda ,\mu }(z)$, defined in terms of Kummer's
confluent hypergeometric functions $M$ and $U$ by
\begin{eqnarray*}
M_{\lambda ,\mu }(z) &=&e^{-\frac{1}{2}z}z^{\frac{1}{2}+\mu }M(\displaystyle%
\frac{1}{2}+\mu -\lambda ,\ 1+2\mu ,\ z), \\
W_{\lambda ,\mu }(z) &=&e^{-\frac{1}{2}z}z^{\frac{1}{2}+\mu }U(\displaystyle%
\frac{1}{2}+\mu -\lambda ,1+2\mu ,\ z).
\end{eqnarray*}

The Whittaker functions $M_{\lambda ,\mu}(z)$ satisfy the following
integral representations

\begin{equation}
M_{\lambda ,\mu }(z)=\displaystyle\frac{\Gamma (1+2\mu )z^{\mu +\frac{1}{2}%
}2^{-2\mu }}{\Gamma (\frac{1}{2}+\mu -\lambda )\Gamma
(\frac{1}{2}+\mu +\lambda )}\int_{-1}^{1}e^{\frac{1}{2}zt}(1+t)^{\mu
-\frac{1}{2}-\lambda
}(1-t)^{\mu -\frac{1}{2}+\lambda }dt,\displaystyle{Re}\mu +\frac{1}{2}%
>|{Re}\lambda |.  \label{Wintrep}
\end{equation}

\section{Definitions and properties of the GOSWFs}


To define the GOSWFs, we 
introduce the following
operator

\begin{equation}  \label{1}
\mathcal{F}_{c}^{(\alpha,\beta)}[\phi](x)=\int_{-1}^{1}e^{c(xy-1)}\phi(y)%
\omega_{\alpha,\beta} (y)dy, \:\: x \in(-1,1), c>0.
\end{equation}
%


Let consider the Hilbert space
$\mathcal{H}=L^2((-1,1),\omega_{\alpha,\beta} (y)dy)$ as the domain
of $\mathcal{F}_{c}^{(\alpha,\beta)}$ and denote by
$$LB_{\omega_{\alpha,\beta}}^c=\left\{f(x)=%
\int_{-1}^{1}e^{c(xy-1)}g(y)\omega_{\alpha,\beta} (y)dy,\:\: g \in
L^2((-1,1),d\omega_{\alpha,
\beta})\right\}=Rg(\mathcal{F}_{c}^{(\alpha,\beta)})$$ its range.
Remark here that straightforward
 computation shows that the operator
 $\mathcal{F}_{c}^{(\alpha,\beta)}$ is a bijection between $L^2((-1,1),\omega_{\alpha,\beta}
(y)dy)$ and $LB_{\omega_{\alpha,\beta}}^c.$

Let now introduce the scalar product defined on the range of
$\mathcal{F}_{c}^{(\alpha,\beta)}$, denoted for simplicity
$Rg(\mathcal{F}_{c}^{(\alpha,\beta)})=\mathcal{\widetilde{H}}$.
Namely, let $f,g\in LB_{\omega_{\alpha,\beta}}^c$, i.e.
$f=\mathcal{F}_{c}^{(\alpha,\beta)}[F]$ and
$g=\mathcal{F}_{c}^{(\alpha,\beta)}[G]$ for some $F,G\in
\mathcal{H}$, then
\[
<f,g>_{\mathcal{\widetilde{H}}}=<F_p,G_p>_{\mathcal{H}},
\]
where  $F_p$ and $G_p$ are the orthogonal projections of $G$ and $F$
to the orthogonal complement of the null-space of
$\mathcal{F}_{c}^{(\alpha,\beta)}$,
$N(\mathcal{F}_{c}^{(\alpha,\beta)})$. Since
$\mathcal{F}_{c}^{(\alpha,\beta)}$ is one to one then
$N(\mathcal{F}_{c}^{(\alpha,\beta)})=\{0\},$ therefore $F=F_p$,
$G=G_p$, and

\begin{equation}\label{scalar-prod}
<f,g>_{\mathcal{\widetilde{H}}}=<F,G>_{\mathcal{H}}.
\end{equation}
Since  $\mathcal{F}_{c}^{(\alpha,\beta)}$ is a self adjoint Hilbert
Shmidt operator, then its eigenfunctions denoted by
$\psi_{n}^{(\alpha,\beta)}(x; c),$ form an orthogonal basis of
$\mathcal{H}\cap
N(\mathcal{F}_{c}^{(\alpha,\beta)})^{\bot}=\mathcal{H}.$

 The GOSWFs
are defined as follows:


\begin{definition}
The GOSWFs are defined as:

\begin{enumerate}
%

\item the eigenfunctions of $\mathcal{Q}_{c}^{(\alpha ,\beta )}=(\mathcal{F}%
_{c}^{(\alpha ,\beta )})^{\ast }\circ \mathcal{F}_{c}^{(\alpha
,\beta )},$ That is
\begin{equation}
\mathcal{Q}_{c}^{(\alpha ,\beta )}[\psi_{n}^{(\alpha ,\beta )}
](x)=\int_{-1}^{1}K_{\alpha ,\beta }(x,y)\psi_{n}^{(\alpha ,\beta )}
(y)\omega _{\alpha ,\beta }(y)dy=\lambda_{n}^{(\alpha ,\beta
)}(c)\psi_{n}^{(\alpha ,\beta )}(x)\:\:\:,x\in (-1,1),c>0. \label{2}
\end{equation}%
Here $\displaystyle{K_{\alpha ,\beta
}(x,y)=e^{-2c}\int_{-1}^{1}e^{ct(x+y)}\omega
_{\alpha ,\beta }(t)dt}\underset{\text{by (\ref{Wintrep})}}{{=}}\displaystyle%
e^{-2c}2^{\alpha +\beta +1}B(\alpha +1,\beta +1)\frac{M_{\frac{\alpha -\beta }{2},%
\frac{\alpha +\beta +1}{2}}(2c(x+y))}{(2c(x+y))^{\frac{\alpha +\beta +2}{2}}}%
$ and $\lambda_{n}^{(\alpha ,\beta )}(c)=(\mu_{n}^{(\alpha ,\beta
)}(c))^2.$
\end{enumerate}
\end{definition}

In the sequel we denote by $\psi_{n}^{(\alpha,\beta)}(x; c)$ the nth
GOSWFs, $-\chi_{n}^{(\alpha,\beta)}(c)$ denotes the nth eigenvalue
of $\mathcal{D}_{x}
$ corresponding to $\psi_{n}^{(\alpha,\beta)}(x; c).$ Let also denote by $%
\mu_{n}^{(\alpha,\beta)}(c)$ the nth eigenvalue of $\mathcal{F}%
_{c}^{(\alpha,\beta)}$ associated to $\psi_{n}^{(\alpha,\beta)}(x;
c).$

 We summarize below some basic properties of the GOSWFs

\begin{proposition}
For any $c>0$ and $\alpha>-1$, $\beta>-1$, we have

\begin{enumerate}
\item[(i)] $\{\psi_{n}^{(\alpha,\beta)}(x; c)\}_{n=0}^{\infty}$ form a complete orthogonal system of $L^2((-1,1),d\omega_{\alpha,%
\beta})$, namely,

\begin{equation}  \label{3}
\int_{-1}^{1}\psi_{n}^{(\alpha,\beta)}(x; c)
\psi_{m}^{(\alpha,\beta)}(x; c)\omega_{\alpha,\beta} (x)dx
=(\mu_{n}^{(\alpha,\beta)}(c))^2\delta_{mn}.
\end{equation}

%
%
%

\item[(v)] $\{\psi_{n}^{(\alpha,\beta)}(x; c)\}_{n=0}^{\infty}$ form a
complete orthogonal system of the bilateral weighted Laplace
bandlimited functions given by,
$$LB_{\omega_{\alpha,\beta}}^c=\left\{f(x)=%
\int_{-1}^{1}e^{c(xy-1)}g(y)\omega_{\alpha,\beta} (y)dy,\:\: g \in
L^2((-1,1),d\omega_{\alpha, \beta})\right\}.$$ 

\item[(vi)] The derivative of $\mu_{n}^{(\alpha,\beta)}(c)$ with respect to $c$ is
given by

\begin{equation}\label{*5454}
\frac{\partial \mu_{n}^{(\alpha,\beta)}(c)}{\partial c} =\frac{1
}{\mu_{n}^{(\alpha,\beta)}(c)}\left(\frac{I_n(c)}{c}-(\mu_{n}^{(\alpha,\beta)}(c))^2\right)
\end{equation}%
where $\displaystyle{I_n(c)=\int_{-1}^{1}v\psi_{n}^{(\alpha,\beta)}(v; c)\frac{\partial \psi_{n}^{(\alpha,\beta)}(v; c)}{%
\partial v}\omega_{\alpha,\beta}(v)dv.}$

\item[(vii)] For all $x\in(-1,1)$ we have the following inequality:
$|\psi_{n}^{(\alpha,\beta)}(x; c)|\leq
\frac{1}{\mathbf{P}_0^{(\alpha,\beta)}(x)}.$
\end{enumerate}
\end{proposition}

\textbf{Proof:} The propertie $(i)$ can be derived from the theory
of Hilbert Shmidt self adjoint operator, see \cite{***}
We restrict ourself here to prove $(v)$ and $(vi).$

Let $f\in LB_{\omega_{\alpha,\beta}}^c,$ then there exists $g \in
L^2((-1,1),d\omega_{\alpha, \beta})$ such that
\begin{equation}\label{*0}
    f(x)=\int_{-1}^{1}e^{c(xy-1)}g(y)\omega_{\alpha,\beta} (y)dy.
\end{equation}
Since  $g \in L^2((-1,1),d\omega_{\alpha, \beta}),$ then by $(i)$ we
have
\begin{equation}\label{*18}
    g(y)=\displaystyle{\sum_{k\in \mathbb{N}}c_k\psi_{k}^{(\alpha,\beta)}(y; c)}.
\end{equation}
Combining (\ref{*0}) and (\ref{*18}) we obtain
\begin{equation}\label{*1}
    f(x)=\displaystyle{\sum_{k\in \mathbb{N}}c_k\mu_{k}^{(\alpha,\beta)}(c)\psi_{k}^{(\alpha,\beta)}(x; c)}.
\end{equation}
This show that the set $\{\psi_{k}^{(\alpha,\beta)}(x; c),\;\;k\in
\mathbb{N}\}$ span the set $LB^c_{\omega_{\alpha,\beta}}.$

to achieve the proof of (v) we remark first that the set
$LB^c_{\omega_{\alpha,\beta}}$ is a reproducing kernel Hilbert space
with .

For the proof of $(vi)$, we adopt the techniques used in
\cite{slepian64} to prove a similar result for the eigenvalue of the
finite Hankel transform. We differentiate both member of the
following equality
\begin{equation}
\int_{-1}^{1}e^{c(vu-1)}\psi_{n}^{(\alpha ,\beta )} (u)\omega
_{\alpha ,\beta }(u)du=\mu_{n}^{(\alpha ,\beta
)}(c)\psi_{n}^{(\alpha ,\beta )}(v)\:\:\:,v\in (-1,1). \label{00*8}
\end{equation}
 with respect to $c $, to obtain
\begin{eqnarray}
&&\int_{-1}^{1}(uv-1)e^{c(uv-1)}\psi_{n}^{(\alpha,\beta)}(u;
c)\omega_{\alpha,\beta}(u)du+\int_{-1}^{1}e^{c(uv-1)}\frac{\partial
\psi_{n}^{(\alpha,\beta)}(u;
c)}{\partial c}%
\omega_{\alpha,\beta}(u)du=  \notag \\
&&\frac{\partial \mu_{n}^{(\alpha,\beta)}(c)}{\partial
c}\psi_{n}^{(\alpha,\beta)}(v;
c)+\mu_{n}^{(\alpha,\beta)}(c)\frac{\partial
\psi_{n}^{(\alpha,\beta)}(v; c)}{\partial c}.  \label{*1}
\end{eqnarray}%
Differentiating (\ref{00*8}) with respect to $v$, one gets
\begin{equation}
\int_{-1}^{1}ue^{c(uv-1)}\psi_{n}^{(\alpha,\beta)}(u;
c)\omega_{\alpha,\beta}(u)du=\frac{\mu_{n}^{(\alpha,\beta)}(c)}{c}\frac{\partial
\psi_{n}^{(\alpha,\beta)}(v; c)}{\partial v}. \label{*2}
\end{equation}%
Combining (\ref{*1}) and (\ref{*2}) to obtain
\begin{eqnarray}
v\frac{\mu_{n}^{(\alpha,\beta)}(c)}{c}\frac{\partial
\psi_{n}^{(\alpha,\beta)}(v; c)}{\partial
v}-\int_{-1}^1e^{c(uv-1)}\psi_{n}^{(\alpha,\beta)}(u;
c)\omega_{\alpha,\beta}(u)du&+&\nonumber\\\int_{-1}^{1}e^{c(uv-1)}\frac{\partial
\psi_{n}^{(\alpha,\beta)}(u;
c)}{%
\partial c}\omega_{\alpha,\beta}(u)du=
\frac{\partial \mu_{n}^{(\alpha,\beta)}(c)}{\partial
c}\psi_{n}^{(\alpha,\beta)}(v;
c)+\mu_{n}^{(\alpha,\beta)}(c)\frac{%
\partial \psi_{n}^{(\alpha,\beta)}(v;
c)}{\partial c}.&&  \label{*3}
\end{eqnarray}%
Multiply both sides of (\ref{*3}) by $\psi _{n}^{(\alpha,\beta)}(v)$
and integrate over $(-1,1) $. One finds
\begin{eqnarray}
&&\frac{\mu_{n}^{(\alpha,\beta)}(c)}{c}\int_{-1}^{1}v\psi_{n}^{(\alpha,\beta)}(v;
c)\frac{\partial \psi_{n}^{(\alpha,\beta)}(v;
c)}{\partial v}%
\omega_{\alpha,\beta}(v)dv-\mu_{n}^{(\alpha,\beta)}(c)\|\psi_{n}^{(\alpha,\beta)}(u;
c)\chi_{(-1,1)}\|_2^2\nonumber\\
&&\int_{-1}^{1}\psi_{n}^{(\alpha,\beta)}(v;
c)\int_{-1}^{1}e^{c(uv-1)}\frac{\partial
\psi_{n}^{(\alpha,\beta)}(u;
c)}{\partial c}\omega_{\alpha,\beta}(u)\omega_{\alpha,\beta}(v)dudv=  \notag \\
&&\frac{\partial \mu_{n}^{(\alpha,\beta)}(c)}{\partial c}\Vert
\psi_{n}^{(\alpha,\beta)}(.; c)\chi _{(-1,1)}\Vert
_{2}^{2}+\mu_{n}^{(\alpha,\beta)}(c)\int_{-1}^{1}\psi_{n}^{(\alpha,\beta)}(v;
c)\frac{\partial \psi_{n}^{(\alpha,\beta)}(v;
c)}{%
\partial c}\omega_{\alpha,\beta}(v)dv.  \label{*4}
\end{eqnarray}%
Using Fubini's Theorem and (\ref{00*8}) together with the
normalization of the GOSWFs, the equality (\ref{*4}) can be simply
written as follows
\begin{equation}
\frac{\partial \mu_{n}^{(\alpha,\beta)}(c)}{\partial c} =\frac{1
}{\mu_{n}^{(\alpha,\beta)}(c)}\left(\frac{I_n(c)}{c}-1\right)
\label{*5}
\end{equation}%
where $\displaystyle{I_n(c)=\int_{-1}^{1}v\psi_{n}^{(\alpha,\beta)}(v; c)\frac{\partial \psi_{n}^{(\alpha,\beta)}(v; c)}{%
\partial v}\omega_{\alpha,\beta}(v)dv.}$

Finally, the proof of (vii) is based on the use of Schwartz
inequality and the GOSWFs normalization.

Remark here that one can check numerically that the eigenvalues
$\mu_{n}^{(\alpha,\beta)}(c)$ decay exponentially to $0.$ This
statement was proved in the special cases of the PSWFs and Wang's
GPSWFs. It will be the subject of a future work.

\section{Numerical computation of the eigenfunctions and the eigenvalues of
the finite bilateral Laplace transform}

To obtain the approximate spectrum of finite bilateral weighted
Laplace transform we use the Gaussian quadrature method. More
precisely, as it is done in \cite{Karoui}, the following theorem
provides a discretization formula for eigenproblem (\ref{00*}) as
well as an interpolation formula for the approximate GOSWFs.

\begin{theorem} Let  $\epsilon $ be a
real number satisfying $0<\epsilon < 1$. Let
$K_\epsilon(\alpha,\beta,c)$ be the integer defined as follows
$$K_\epsilon(\alpha,\beta,c)=\inf\{K \in\mathbb{N},
\frac{e^{2c}}{\sqrt{2K\pi}}\frac{2^{\alpha+\beta+1}
K!\Gamma(K+\alpha+\beta+1)\Gamma(K+\alpha+1)\Gamma(K+\beta+1)}{(2K+\alpha+\beta+1)(\Gamma(2K+\alpha+\beta+1))^2}\leq\epsilon
|\mu_n^{(\alpha,\beta)}|\}.$$ Then for any
$k\geq\max([2ec]+1,K_{\epsilon }),$ we have

\begin{equation}\label{numres}
\displaystyle\sup_{x\in [ 0,1]}\left|\psi_n^{(\alpha,\beta)}
(x;c)-\frac{1}{\mu_n^{(\alpha,\beta)}(c)}\sum_{j=1}^{k}\omega_{j}
e^{(cxy_{j})}\psi_n^{(\alpha,\beta)} (y_{j};c)\right|<\epsilon
\end{equation}

Here, $(y_{j})_{1\leq j\leq n}$, denote the different zeros of
$\mathbf{P}_{n}^{\alpha ,\beta }(x)$.
\end{theorem}

 \begin{remark}
Numerical evidences show that $\mu _{n}^{(\alpha,\beta)}(c)$ decay
rapidly to  $0$ when $n$ goes to $+\infty.$
 \end{remark}

\section{Applications of the GOSWFs and numerical results}
In this section, we give two applications of the GOSWFs: The fisrt
application, deals with the approximation of bilateral Laplace
band-limited  by the use of the GOSWFs. As a second application of
the GOSWFs, we use them to invert the bilateral Laplace transform of
time limited functions. We provide the reader by several numerical
results that illustrate these applications as well as of the
eigenvalues $\mu _{n}^{(\alpha,\beta)}(c)$ for different values of
$\alpha,$ $\beta$ and $c$  and some curves of GOSWFs .

\subsection{GOSWFs and approximation of bilateral
Laplace band-limited signals}

{\color{blue} Similarly to was done in \cite{Slepian-Pollak}, for
the PSWFs, and in
 \cite{Moumni1} for the CPSWFs, we can show that the GOSWFs solve the problem of
signal concentration energy, which is an important problem from
signal processing. More precisely, we assert that among the set
$LB_{\omega_{\alpha,\beta}}^c$ of bilateral Laplace band-limited
signals with bandwidth $c>0,$ $\psi_{0}^{(\alpha,\beta)}(.;c)$ has
the most concentrated energy on $(-1,1)$. More generally, for any
integer $n\geq 1,\: \psi_{n}^{(\alpha,\beta)}(.;c)$ is the $(n+1)$th
most concentrated signal in $(-1,1)$ which is orthogonal to $
\psi_{0}^{(\alpha,\beta)}(.;c),\ \ldots,\:
\psi_{n}^{(\alpha,\beta)}(.;c) $. By the properties of the GOSWFs
given in Proposition~1, for any  $f\in
L^{2}((-1,1),d\omega_{\alpha,\beta})$, the GOSWFs-based expansion
formula of $f$ over $(-1,1)$ is given by $f(t)=\displaystyle
\sum_{n\geq
0}\alpha_{n}\psi_{n}^{(\alpha,\beta)}(t;c)$. since $\displaystyle \int_{-1}^{1}\psi_{n}^{(\alpha,\beta)}(t;c)%
\psi_{m}^{(\alpha,\beta)}(t;c)\omega_{\alpha,\beta}(t)dt=(\mu_n^{(\alpha,\beta)}(c))^2\delta_{nm}$,
then

$$\displaystyle \left|
\alpha_{n}\right|=\frac{1}{(\mu_n^{(\alpha,\beta)}(c))^2}\left|\int_{-1}^{1}f(t)\psi_{n}^{(\alpha,\beta)}(t;c)\omega_{\alpha,\beta}(t)dt\right|\leq\frac{\|f\|_2}{|\mu_n^{(\alpha,\beta)}(c)|}$$.

Let $f_{N}^{GOSWFs}$ denotes the N-term truncated GOSWFs series
expansion of $f$, given by

$f_{N}^{GOSWFs}(t)=\displaystyle
\sum_{n=0}^{N}\alpha_{n}\psi_{n}^{(\alpha,\beta)}(t;c),\
t\in(-1,1)$.

Since the sequence $\mu_n^{(\alpha,\beta)}(c)$ decay rapidly to $0$ when $n$ goes to $\infty$, and since $%
\Vert\psi_{n,c}\chi_{(-1,1)}\Vert_{\infty}\leq\frac{1}{\mathbf{P}_0^{(\alpha,\beta)}(x)}$,
then $f_{N}^{GOSWFs}$ converges rapidly to $f$ and
$|f_{N}^{GOSWFs}(t)-f(t)|=O(\mu_n^{(\alpha,\beta)}(c)),\ \forall
t\in(-1,1)$. Moreover, if $f\in
L^{2}((-1,1),d\omega_{\alpha,\beta})$, then
$\displaystyle\lim_{N\rightarrow+ \infty}\Vert
f-f_{N}^{GOSWFs}\Vert_{2}=0$. In the sequel we denote by
$f_{N}^{Jacobi}$  the N-term truncated Jacobi series expansion of
$f$, given by

$f_{N}^{Jacobi}(t)=\displaystyle
\sum_{n=0}^{N}\alpha_{n}\mathbf{P}_n^{(\alpha,\beta)}(t;c),\
t\in(-1,1)$. We compare numerically, for a given $N$ and a given
function $f\in {LB^c_\omega}_{(\alpha,\beta)},$ which function among
$f_{N}^{Jacobi}$ or $f_{N}^{GOSWFs}$ approach the best $f.$

}

\subsection{Numerical results}

In this subsection we give several numerical results that illustrate
the theoretical results of the previous sections. We have considered
different values of the bandwidth $c$ and  the parameters $\alpha$
and $\beta.$ Also, we have applied the Gaussian quadrature based
method for the computation of the spectrum and the eigenfunctions of
the finite weighted bilateral  Laplace transform
$\mathcal{F}_{c}^{(\alpha,\beta)}$ with $N=40 $ quadrature points.
In table 1 we have listed  the obtained eigenvalues $\mu
_{n}^{(\alpha,\beta)}(c)$ with different values of the parameter
$c,$ $\alpha$ and $\beta.$

\begin{center}
\begin{table}[h!]
\begin{tabular}{ccccc}
\hline n & $c=1$                 & $c=1$                 & $c=6$                & $c=6$ \\
\hline & $(\alpha,\beta)=(0,0)$  & $(\alpha,\beta)=(3,3)$  & $(\alpha,\beta)=(6,7)$  & $(\alpha,\beta)=(5,5)$ \\
\hline\\
0 & $0.779836289$              & $0.338455158        $      &$ 0.199353974*10^{-2}$     &$ 0.211037689*10^{-2}$ \\
2 & $0.328060086*10^{-1}$       & $0.305509085*10^{-2} $     &$ 0.242164507*10^{-3}$     &$ 0.409615392*10^{-3}$ \\
4 & $0.178076210*10^{-3}$       &$ 0.982704157*10^{-5}  $    &$ 0.146846578*10^{-4}$     &$ 0.312721188*10^{-4}$\\
6 & $3.771944953*10^{-7}$       & $1.571133391*10^{-8}   $   &$ 5.322716251*10^{-7}$     &$ 0.134135193*10^{-5}$ \\
8 & $4.247451396*10^{-10}$      &$ 1.481487170*10^{-11}   $  &$ 1.279499739*10^{-8}$     &$ 3.673872446*10^{-8}$  \\
10 &$2.966038648*10^{-13}$      & $9.151623350*10^{-15}    $ &$ 2.179482513*10^{-10}$    &$ 6.946790778*10^{-10}$\\
15 &$8.099510876*10^{-22}$      &$ 2.071915188*10^{-23}$     &$ 2.426475139*10^{-15}$    &$ 9.333247151*10^{-15}$ \\
20 &$4.268133206*10^{-31}  $    & $9.817546181*10^{-33}$     &$ 6.690930862*10^{-21}$    &$ 2.916259632*10^{-20}$  \\
\hline
\end{tabular}
\vskip 0.2cm \hspace{3.5cm} \caption{Values of the eigenvalues
$\protect\mu_{n}^{(\alpha,\beta)}(c)$ of the finite bilateral
weighted Laplace transform $\mathcal{F}_{c}^{(\alpha,\beta)}$
corresponding to  different values of parameter $\alpha,$ $\beta$
and $c.$}
\end{table}
\end{center}

We remark easily that the GOSWFs are reduced to the PSWFs in the
special case where the integral operator
$\mathcal{F}_c^{(\alpha,\beta)}$ is replaced by
$\mathcal{\widetilde{F}}_c=e^{ic}\mathcal{F}_{ic}^{(0,0)},$ Here
$i^2=-1.$
 Hence the eigenvalues of
the GOSWFs coincides with those of the PSWFs. In  table~2, we give
some results of the eigenvalues of
$\frac{c}{2\pi}\mathcal{\widetilde{F}}_c\circ(\mathcal{\widetilde{F}}_c)^*,$
such values was given also in \cite{karoui-moumni1}.
\begin{center}
\begin{table}[h!]
\begin{tabular}{cccc}
\hline n & $\widetilde{c}=2$   & $\widetilde{c}=4$       &
$\widetilde{c}=6$\\
\hline \\
0  & $0.8805599223*10^{0}$              & $0.9958854904*10^{0}$  &$0.9999018826*10^{0}$\\
 5  &$1.9358522020*10^{-7}$       & $0.3812917217*10^{-3}$        &$0.2738716624*10^{-1}$\\
  10 &$2.1680118965*10^{-19}$      &$4.5252284693*10^{-13}$       &$2.2189805452*10^{-9}$\\
   15 &$1.6563615010*10^{-33}$      &$3.5519079602*10^{24}$      &$1.0163838373*10^{-18}$ \\
   20 &$4.7105458228*10^{-49}$      &$1.0352225590*10^{-36}$   &$1.7132439301*10^{-29S}$\\
   \hline
\end{tabular}
\vskip 0.2cm \hspace{3.5cm} \caption{Values of the eigenvalues of
$\frac{c}{2\pi}\mathcal{\widetilde{F}}_c\circ(\mathcal{\widetilde{F}}_c)^*,$
corresponding to the different values of parameter $c.$}
\end{table}
\end{center}

For the computation of the GOSWFs over $(-1,1)$ by the first method,
we have used the equality (\ref{18}) with a maximum truncation order
$K=10,$  see Figure~1. For the computation of the GOSWFs over
$(-1,1)$ by the second method, we have used the equality
(\ref{numres}) with a maximum truncation order $K=10$ quadrature
points, see Figure~2 and Figure~3. Remark, that for $\alpha=\beta,$
$\psi_n^{(\alpha,\beta)}(.;c)$ has the same parity as its order $n$
contrarily to the case where $\alpha\neq\beta.$

\begin{figure}[h!]
\begin{center}
\includegraphics[width=15cm,height=7cm]{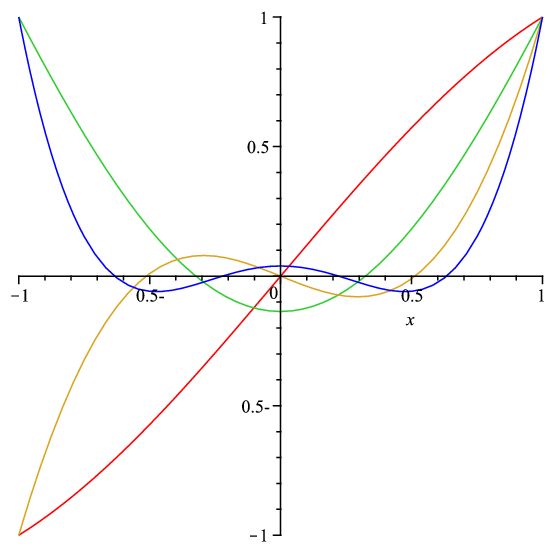}
\end{center}
\vspace{1.5cm} \caption{Graphs of $\psi_{k}^{(\alpha,\beta)}(.;c),$
$k=1,2,3,4$ associated to the parameter $c=2$ and $\alpha=\beta=3$
by the differential operator based method} \label{fig1}
\end{figure}

\vspace{0.5cm}
\begin{figure}[h!]
\begin{center}
\includegraphics[width=15cm,height=7cm]{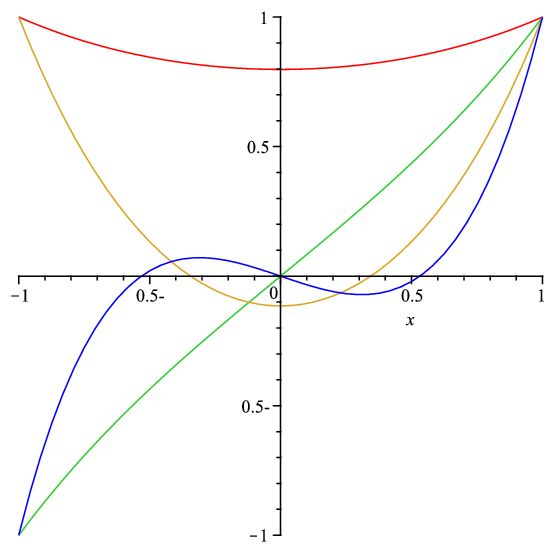}
\end{center}
\vspace{1.5cm} \caption{Graphs of the first four GOSWFs associated
to the parameter $c=2$ and $\alpha=3$ and $\beta=3$ by the Gaussian
quadrature based method} \label{fig2}
\end{figure}

\vspace{0.5cm}
\begin{figure}[h!]
\begin{center}
\includegraphics[width=15cm,height=7cm]{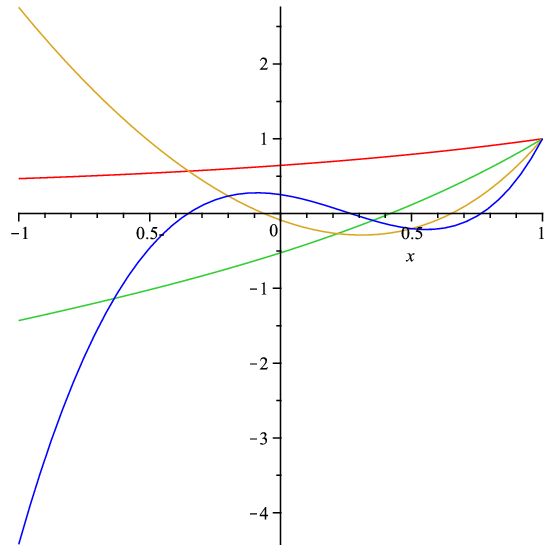}
\end{center}
\vspace{1.5cm} \caption{Graphs of the first four GOSWFs associated
to the parameter $c=1$ and $\alpha=2$ and $\beta=3$ by the Gaussian
quadrature based method} \label{fig3}
\end{figure}

As a comparison  between the approximation of functions in
$LB_\omega$ by the use of the GOSWFs and Jacobi polynomials, we give
the curves of three functions in $LB_\omega,$ their approximation by
the GOSWFs and also by the Jacobi polynomials, see figures~4-6.

\begin{figure}[h!]
   \includegraphics[width=15cm,height=8cm]{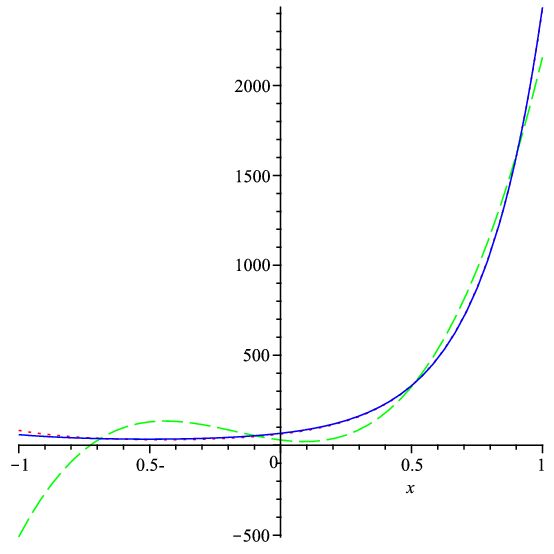}%
   \caption{Curves of $F(x)=\frac{e^{-cx}}{c}\left(\frac{L!}{(-x)^{L+1}}-e^{2cx}\sum_{j=0}^L\frac{L!}{j!}\frac{(2c)^j}{(-x)^{L-j+1}}\right)$
 (line) and its
approximated signal $F_N^{GOSWFs} $ (dot), $F_N^{Jacobi} $ (dash) ,
where $N=3,$ $c=5,$ $\alpha=0,$ $\beta=1,$ and $L=2.$}
\end{figure}

\begin{figure}[h!]
   \includegraphics[width=15cm,height=8cm]{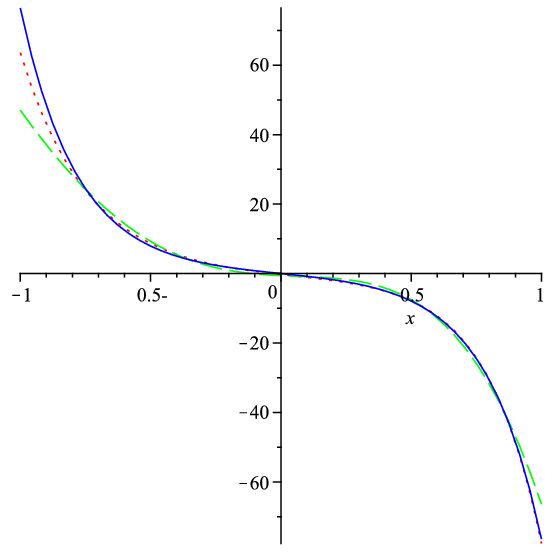}%
   \caption{Curves of $F(x)=\pi I_1(-cx)$
 (line) and its
approximated signal $F_N^{GOSWFs} $ (dot), $F_N^{Jacobi} $ (dash) ,
where $N=4,$ $c=5,$ $\alpha=1,$ $\beta=2.$}
\end{figure}

\begin{figure}[h!]
   \includegraphics[width=15cm,height=8cm]{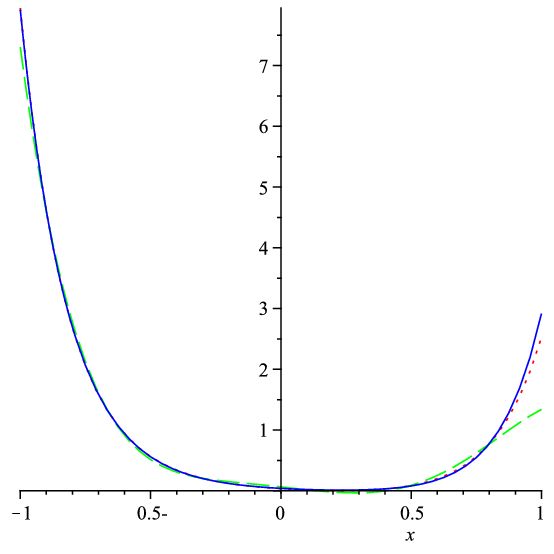}%
   \caption{Curves of $F(x)=\frac{e^{-cx}}{c}\left((-1)^{\nu}*\sqrt{0.5\pi x}I_{\nu}(cx)K_{\nu}(\nu+0.5,cx)\right)$
 (line) and its
approximated signal $F_N^{GOSWFs} $ (dot), $F_N^{Jacobi} $ (dash) ,
where $N=5,$ $\nu=2,$ $c=6,$ $\alpha=2,$ $\beta=1.$}
\end{figure}

\end{document}